\documentclass{amsart}
\usepackage[T1]{fontenc}
\begin{document}

\centerline{GEOMETRIC PROPERTIES OF SOME BANACH ALGEBRAS RELATED TO THE}

\medskip

\centerline{FOURIER ALGEBRA ON LOCALLY COMPACT GROUPS.}

\medskip

\centerline{EDMOND E. GRANIRER}

\bigskip

\noindent ABSTRACT. {\it{Let}} $A_p(G)$ {\it{denote the Figa-Talamanca-Herz Banach Algebra of the locally compact group}} $G$, {\it{thus}} $A_2(G)$ {\it{is the Fourier Algebra of}} $G$. {\it{If}} $G$ {\it{is commutative then}} $A_2(G)=L^1(\hat{G})^\wedge$. {\it{Let}} $A_p^r(G)=A_p\cap L^r(G)$ {\it{with norm}} $\| u\|_{A_p^r}=\| u\|_{A_p}+\| u\|_{L^r}$.

\noindent {\it{We investigate for which}} $p$, $r$, {\it{and}} $G$ {\it{do the Banach algebras}} $A_p^r(G)$ {\it{have the Banach space geometric properties: The Radon-Nikodym Property (RNP), the Schur Property (SP) or the Dunford-Pettis Property (DPP).}}

{\it{The results are new even if}} $G=R$ {\it{(the real line) or}} $G=Z$ {\it{(the additive integers)}}.

\medskip

\noindent INTRODUCTION. Let $G$ be a locally compact group and let $A_p(G)$ denote the Figa-Talamanca-Herz Banach algebra of $G$, as defined in \cite{Her1}, thus generated by $L^{p'}\ast L^{\vee p}(G)$, where $1<p<\infty$, and $1/p+1/p'=1$, see sequel. Hence $A_2(G)$ is the Fourier algebra of $G$ as defined and studied by Eymard in \cite{Eym}. If $G$ is abelian then $A_2(G)=L^1(\hat{G})^\wedge$.  

Denote $A_p^r(G)=A_p\cap L^r(G)$, for $1\leq r\leq\infty$, equipped with the norm
\begin{equation}  
 \| u\|_{A_p^r}=\| u\|_{A_p}+\| u\|_{L^r}.\qquad\mbox{If $r=\infty$ let $A_p^r(G)=A_p(G)$.}
\end{equation}
						
If $G$ is abelian then, $A_2^r(G)=L^1(\hat{G})^\wedge \cap L^r(G)$ with the norm
\begin{equation}
\| u\| =\| f\|_{L^1(\hat{G})}+\|\hat{f}\|_{L^r(G)}\quad\mbox{if }u=\hat{f}.
\end{equation}

The study of these Banach Algebras started in a beautiful paper of Larsen Liu and Wang \cite{Lar} in the abelian case, and continued in \cite{Lai}, \cite{Lai2}, \cite{Gr1}-\cite{Gr5} ... etc.

Let $X$ be a Banach space. Then 

\noindent $X$ {\it{has the Schur Property (SP) if weak convergent sequences are norm convergent}}.

\noindent $X$ {\it{has the Dunford-Pettis Property (DPP) if whenever}} $(x_n),(x_n^\ast )$ {\it{are weakly null sequences in}} $X$ {\it{and}} $X^\ast$ {\it{respectively, then}}, $\lim (x_n,x_n^\ast )=0$.

{\it{Clearly the SP implies the DPP,}} \cite{Dtl}.

The Banach space $\ell^1$ is a dual Banach space which has the SP, hence the DPP. But for any measure space, $L^1(\mu )$, {\it{has the}} DPP, \cite{Dtl} p.19, yet it {\it{does not have}} the SP if $\mu$ is non atomic.

$X$ {\it{has the Krein-Milman (KMP) property if any closed convex bounded subset}} $B$ {\it{is the norm closed convex hull of its extreme points (ext B)}}

$X$ {\it{has the Radon Nikodym property (RNP) if every such}} $B$ {\it{is the norm closed convex hull of its strongly exposed points (strexp B), see sequel or \cite{Die} p.190 and p.218}}. 

\renewcommand{\thefootnote}{}

\newpage 

Points in {\it{strexp(B)}} are points of {\it{ext(B)}} that have beautiful smoothness properties. In particular they are weak to norm continuity points of $B$ and are peak points of $B$, see sequel.

\noindent Quoting Jerry Uhl: ``{\it{A Banach space has the RNP if it's unit ball ``wants to be weakly compact, but just cannot make it''}}.

\noindent Denote by $PM_p(G)=A_p(G)^\ast$, and by $PF_p(G)$, the norm closure in $PM_p(G)$ of $L^1(G)$, (as a space of left convolutors on $L_p(G)$). Let $W_p(G)=PF_p(G)^\ast$

Then $W_p(G)$ is a Banach algebra of bounded continuous functions on $G$ containing the ideal $A_p(G)$, see Cowling \cite{Cow}.

Let $W_p^r(G)=W_p\cap L^r(G)$, with the sum norm.

If $G$ is abelian and $p=2$ then $W_2(G)=M(\hat{G})^\wedge $, where $M(G)$ is the space of bounded Borel measures on $G$.

\bigskip

\noindent {\it{Our Main Result on the SP and the DPP, is the}}

\medskip

\noindent {\it{THEOREM A:}} {\it{Let G be a noncompact locally compact group}}.
\begin{enumerate}
\item $\forall 1< r\leq\infty\ ,A_p^r(G)$ {\it{does not have the SP}}.

\noindent {\it{Yet}} $A_p^1(G)$ {\it{has the SP and the RNP, if G is discrete.}}

\item {\it{If G is weakly amenable then}} $\forall 1<r\leq p^\prime$, $A_p^r(G)$ {\it{does not have even the DPP yet it has the RNP}}.
\end{enumerate}

The above result is of interest even if $G=R$ or $G=Z$.

\medskip

\noindent {\it{CONJECTURE:}} {\it{If G is a connected semisimple Lie group with finite center then}} $A_2^r(G)$ {\it{has the DPP for any}} $r>2$.

\medskip

\noindent QUESTION: Does this hold if $G=Z$ or $G=R$?

It has been proved by W.~Braun, in an unpublished preprint \cite{Bra}, {\it{that if G is amenable, then}} $A_p^1(G)$ {\it{is a dual Banach space with the RNP}}. 

{\it{This is improved in our Main Result on the RNP}} in the:

\medskip

\noindent {\it{THEOREM B:}} {\it{Let G be a weakly amenable (see sequel) locally compact group. Then}}
\begin{enumerate}
\item $\forall 1\leq r\leq p^\prime ,A_p^r(G)=W_p^r(G) \mbox{ and }$\hfill\break

\hspace{0.05mm} $A_p^r(G)$ {\it{is a dual Banach algebra with the RNP}}. 

\noindent {\it{If G is unimodular, this is the case}} $\forall 1\leq r\leq\max (p,p^\prime )$.

\item {\it{If G is a noncompact connected semisimple Lie group with finite center then}} $\forall r>2$, $A_2^r(G)$, {\it{does not have the RNP and is not a dual Banach space.}}
\end{enumerate}

The above is new and of interest even if $G=R$ or $G=Z$.

\medskip

\noindent QUESTION: Does (2) hold true if $G=Z$, or $G=R$?

\medskip

\noindent {\it{THEOREM C:}} {\it{Assume that G is second countable and weakly amenable, and}} $1<p<\infty$. 

{\it{If for some}} $t\leq\infty$, $A_p^t(G)$ has the RNP, then so does $A_p^r(G)$, $\forall 1\leq r\leq t$

Fell groups, see sequel or \cite{Bag}, are non compact groups for which $A_2(G)$ has the RNP. Hence $A_p^t(G)$ has the RNP for all $1\leq t\leq\infty$, for such $G$.

\medskip

\noindent This paper was inspired by the important paper of F.~Lust-Piquard \cite{Lus}, where the RNP and SP were investigated for $PM_p(E)$ for nowhere dense compact sets $E$, for abelian $G$. See also T.~Miao and P.F.~Mah \cite{Mah}.

\bigskip

\noindent NOTATIONS AND DEFINITIONS: Denote as in \cite{Her1}
\begin{equation}
A_p(G)=\left\{ u=\Sigma u_n\ast v_n^\vee ;
   u_n\in L^{p^\prime}, v_n\in L^p,
		  \Sigma \| u_n\|_{L^{p^\prime}}
			\| v_n\|_{L^p}<\infty\right\} ,
\end{equation}
where the norm of $u\in A_p$ is the infimum of the last sum over all the representations of $u$ as above.

Let $C_0(G)[C_c(G)]$ denote the continuous functions which tend to $0$ at $\infty$, [with compact support], with norm $\| u\|_\infty =\sup\left\{ |u(x)|; x\in G\right\}$.

The {\it{group G is weakly amenable}} if $A_2(G)$ has an approximate identity $\{ v_\alpha\}$ bounded in the norm of $B_2(G)$, the space of Herz-Schur multipliers, see, \cite{Eyma} (or \cite{deC}, \cite{Gr5}).

\noindent Any closed subgroup $G$ of any finite extension of the general Lorenz group $SO_0(n,1)$, for $n>1$, hence the free group on $n>1$ generators is weakly amenable but non amenable. For this and much more see \cite{deC}.

\medskip

\indent {\it{Definition:}} {\it{Let B be a bounded subset of the Banach space X and}} $b\in B$. $b$ {\it{is a strongly exposed point of B (and strexp(B) denotes the set of all such), if}} $\exists b^\ast\in X^\ast$ {\it{such that:}} $\mbox{Re}b^\ast (x)<\mbox{Re}b^\ast (b)$, $\forall x\in B${\it{ and }} $x\not= b$,{\it{ and}}
\begin{equation}
\mbox{Re}b^\ast (x_n)\rightarrow\mbox{Re}b^\ast (b)\mbox{{\it{ for }}}x_n\in B\quad implies\quad 
   \| x_n-b\|\rightarrow 0.\quad\mbox{{\it{(see \cite{Die} p.138)}}}
\end{equation}

{\it{Hence in order to Test an Algorithm for some b in strexp(B) it is Enough to Test it on One Particular Element of}} $X^\ast$.

\bigskip

\centerline{2. MAIN RESULTS}

\noindent (I) THE RNP CASE.

We first improve results in \cite{Gr4}, \cite{Gr5}, by removing the unimodularity of the group in the next

\medskip

\noindent {\it{PROPOSITION 1:}} {\it{Let G be a locally compact group. If}} $p=2$, {\it{or if G is weakly amenable and}} $1<p<\infty$, {\it{then}}
\begin{enumerate}
\item  $(^\ast )W_p\cap L^r(G)=A_p\cap L^r(G),\quad \forall 1\leq r\leq p^\prime$.

{\it{and}} $A_p^r(G)$ {\it{is a dual Banach space.}}

{\it{If G is unimodular then this holds for}} $\forall 1\leq r\leq \max (p,p^\prime )$.

\item {\it{If}} $G=SL(2,R)$ {\it{and}} $p=2$ {\it{then}} $(^\ast )$ {\it{does not hold for any}} $r>2$, {\it{and}} $A_2^r(G)$ {\it{is not a dual space for}} $r>2$, {\it{(see Prop.3)}}.
\end{enumerate}

\medskip

\noindent {\it{REMARK:}} {\it{The interval}} $[1,p^\prime ]$ {\it{is the best one can do even for G=Z and p=2 as proved in \cite{Hew}, (see \cite{Gr4} p.4379, or \cite{Roo})}}.

\medskip

\noindent {\it{PROOF:}} By weak amenability, for all $1<p<\infty$, the $W_p$ norm restricted to $A_p$ is equivalent to the $A_p$ norm, (\cite{Gr5} Corollary 3.7.). If $p=2$ then Kaplansky's density theorem will yield the same result.

Now with the notations of \cite{Gr4} Thm.~2.1. p.4379, if $e_\alpha\in C_c(G)$ is an approximate identity for $L_1(G)$, such that each $e_\alpha$ is the ``square of a special operator'', a la Fendler \cite{Fen} p.129, we have, by the Lemma, loc. cit. that
\begin{equation}
\| e_\alpha\ast w-w\|_{W_{p^\prime}}\rightarrow 0\quad 
   \forall w\in W_{p^\prime}\quad\mbox{afortiori}\quad
	    \forall w\in W_{p^\prime}\cap L^{p^\prime\vee}.
\end{equation}
But, since $e_\alpha\in C_c(G)$, we have for such $w$, that
\begin{equation}
e_\alpha\ast w\in L^p\ast L^{p^\prime\vee}\subset A_{p^\prime},
   \mbox{ thus $e_\alpha \ast w$ is a Cauchy sequence in $A_{p^\prime}$}.
\end{equation}
Hence $w\in A_{p^\prime}$. It follows that $W_{p^\prime}\cap L^{p^\prime \vee}=A_{p^\prime}\cap L^{p^\prime \vee}$.

However by \cite{Cow} p.91, $W_p=W_{p^{\prime\vee}}$, $A_p=A_{p^{\prime\vee}}$. Hence
\begin{equation}
W_p\cap L^{p^\prime}=A_p\cap L^{p^\prime},\quad\forall 1<p<\infty .
\end{equation}

But $W_p$ contains only bounded functions, hence
\begin{equation}
\forall r\leq p^\prime , W_p\cap L^r=W_p\cap L^{p^\prime}\cap L^r=A_p\cap L^{p^\prime}\cap L^r=A_p\cap L^r.\mbox{ Thus}
\end{equation}
\begin{equation}
\mbox{(i)}\quad W_p\cap L^r=A_p\cap L^r,\quad\forall r\leq p^\prime .
\end{equation}
If $G$ is unimodular then, since $\left( W_p\cap L^r\right)^\vee =\left( A_p \cap L^r\right)^\vee$, it follows that $W_{p^\prime}\cap L^r=A_{p^\prime}\cap L^r$, $\forall r\leq p^\prime$, which holds for all $1<p^\prime <\infty$.

Replace now $p^\prime$ by $p$, then $W_p\cap L^r=A_p\cap L^r$, $\forall r\leq p$.

The above implies the unimodular case.

\noindent By Theorem~2{.}2 of \cite{Gr5} $W_p(G)\cap L^r(G)$ is a dual Banach space for all $1<p<\infty$ and $1\leq r\leq\infty$, and all locally compact groups $G$. This proves~($^\ast$). The proof of (ii) and Theorem~B requires the next results. $\quad\square$

\medskip

\noindent {\it{LEMMA 2:}} {\it{Let}} $G$ {\it{be a locally compact group. Assume that}} $A_p(G)$ {\it{has an approximate identity}} $u_\alpha$ {\it{such that}} $\sup\| u_\alpha\|_\infty\leq B<\infty$. {\it{Then}}
\begin{enumerate}
\item $A_p\cap C_c$ {\it{is norm dense in}} $A_p^r$ {\it{and}}

\item {\it{If G is second countable then}} $A_p^r$ {\it{is norm separable.}}
\end{enumerate}

\medskip

\noindent {\it{PROOF:}}  (1) Let $e_\alpha\in A_p\cap C_c$ satisfy $\| e_\alpha -u_\alpha\|_{A_p}\rightarrow 0$ and 
\begin{equation}
\| e_\alpha -u_\alpha\|_{A_p}\leq 1,\, \forall\alpha .\mbox{ Then }
   \| e_\alpha\|_\infty\leq\| e_\alpha -u_\alpha\|_\infty 
	    +\| u_\alpha\|_\infty\leq 1+B.\mbox{ Hence}
\end{equation}
\begin{equation}
\| e_\alpha v-v\|_{A_p}\leq\| (e_\alpha -u_\alpha )v\|_{A_p}
   +\| u_\alpha v-v\|_{A_p}\rightarrow 0,\, \forall v\in A_p.
	   \mbox{ But if }w\in A_p^r
\end{equation}
and $K\subset G$ is compact such that $\int\limits_{G\sim K}|w|^r\, dx<\epsilon$ then
\begin{equation}
\int\limits_{G\sim K}\left| (e_\alpha -1)w\right|^r
  \leq\int\limits_{G\sim K}(2+B)|w|^r\leq (2+B)\epsilon .
		 \mbox{ But }\int\limits_K\left| (e_\alpha -1)w\right|^r\rightarrow\infty .   
\end{equation}
It thus follows that $\| e_\alpha w-w\|_{A_p^r}\rightarrow 0$. But $e_\alpha w\in A_p\cap C_c$.

\indent (2) $A_p(G)$ is norm separable, hence so is $A_p[K]=\{ u\in A_p(G); spt\, u\subset K\}$, where $K\subset G$. Let $A_p^r[K]=\{ u\in A_p^r(G); spt\, u\subset K\}$. If $K$ is compact then the identity $I:A_p^r[K]\rightarrow A_p[K]$ is 1-1, onto and continuous, hence it is bicontinuous. Hence $A_p^r[K]$ is separable. Let now $K_n\subset\mbox{ int }K_{n+1}\subset G$, be compact (int denotes interior), such that $\cup K_n=G$. It is hence enough to show that $\cup A_p^r[K_n]$ is norm dense in $A_P^r(G)$. 

\noindent By (a) we know that $A_p\cap C_c$ is norm dense in $A_p^r(G)$. But if $v\in A_p^r(G)$ has compact support $S$ then $S\subset K_j$ for some $j$, hence $v\in A_p^r[K_j]$. Thus $\cup A_p^r[K_n]$ is norm dense in $A_p^r(G)$. $\quad\square$

\medskip

\noindent {\it{REMARK:}} We do not know if, $A_p\cap C_c(G)$ is norm dense in $A_p^r(G)$ even for $G=SL(2,R)\triangleleft R^2$, if $p=2$ and any $r$. As shown in \cite{Dor}, $A_2(G)$ does not have an approximate identity, bounded in the multiplier norm.

\medskip

\noindent {\it{COROLLARY 3:}} {\it{Let}} $G$ {\it{be a second countable locally compact group.}} {\it{If}} $G$ {\it{is weakly amenable then}} $\forall 1\leq r\leq p^\prime$, $A_p^r(G)$ {\it{is a separable dual Banach algebra and thus has the RNP.}}

\indent {\it{If G is unimodular, this is the case for}} $1\leq r\leq\max (p,p^\prime )$.

\medskip

\noindent {\it{REMARK:}} {\it{Weak amenability, namely the existence in}} $A_2$ {\it{of an approximate identity norm bounded in}} $B_2$  {\it{depends only on}} $p=2$, {\it{yet the result holds for all}} $p$. {\it{Since by Furuta’s Thm.~2.4 in \cite{Fur},}} $B_2\subset B_p$ {\it{contractively, see also \cite{Gr5} p.23. The}} $B_p$ {\it{norm dominates the multiplier norm by \cite{Fur}.}}

\medskip

\noindent {\it{PROOF:}} $A_p^r(G)$ is a dual Banach space $\forall 1\leq r\leq p^\prime$. But since $G$ is weakly amenable $\forall 1< p<\infty$, $A_p(G)$, has a multiplier norm, bounded approximate identity, by the Remark above. It thus follows by the Lemma above, that $A_p^r(G)$ is norm separable.  But separable dual Banach spaces have the RNP by \cite{Die} p.218. $\quad\square$

\medskip
     
\indent The second countability of $G$ is removed in the main result of this section, namely

\medskip

\noindent {\it{THEOREM B:}} {\it{Let G be a weakly amenable locally compact group and}} $1<p<\infty$. {\it{Then}}
\begin{enumerate}
\item $\forall 1\leq r\leq p^\prime$, $A_p^r(G)=W_p^r(G) \mbox{ and}$

\hspace{0.5mm} $A_p^r(G)$ {\it{is a dual Banach algebra with the RNP.}}
				
\indent {\it{If G is unimodular, this is the case}} $\forall 1\leq r\leq\max (p,p^\prime )$.

\item  {\it{If G is a noncompact connected semisimple Lie group with finite center then $\forall r>2$, $A_2^r(G)$ does not have the RNP and is not a dual Banach algebra.}}
\end{enumerate}

\medskip

\noindent {\it{PROOF:}} (1) By \cite{Die} it is enough to prove that every separable subspace of $A_p^r(G)$ has the RNP. Based on the above Corollary follow the proof of Theorem 3{.}1 on p.22-24 of \cite{Gr5} and \cite{Gr4} p.4381.
   
\indent (2)  By a deep result of Cowling \cite{Clg}, $A_2^r(G)=A_2(G)$ if $r>2$.

\noindent Assume that $A_2(G)$ has the RNP. Then the regular representation is the direct sum of irreducible unitary representations, by K.~Taylor’s \cite{Tlr} Thm.~4.1. Denote by $\hat{G}_r$ the set of all such. Then by \cite{Dix} 14.1.2, 14.3.2, $\hat{G}_r$ contains only square integrable representations. Now, by Lipsman, \cite{Lip} p.412--413, $\hat{G}$ induces the discrete topology on $\hat{G}_r$ (this being the set of all square integrable representations). But by the Corollary on p.228 of \cite{Fel}, the topology of $\hat{G}$ is second countable, since $G$ is such. Hence so is that of $\hat{G}_r$, which is in addition discrete and thus is countable. But by Baggett’s \cite{Bag} Prop. 2.2, a connected semisimple Lie group whose reduced dual is countable is compact.
     
$A_2(G)$ is separable. If it was a dual space it would have the RNP, see \cite{Die}.         

We note that (2) has been proved for $SL(2,R)$ in \cite{Gr4}, by using the support of its Plancherell measure. $\quad\square$

\medskip

\noindent {\it{REMARK:}} {\it{Any closed subgroup of any finite extension of the general Lorenz group}} $SO_0(n,1)$ {\it{for}} $n>1$, {\it{hence the free group of}} $n>1$ {\it{generators (a nonamenable group), is a weakly amenable group.}}

\noindent {\it{This group is a noncompact connected simple Lie group, see}} \cite{deC} {\it{p.474 for this and much more.}}

\bigskip

\noindent (II) INTERVALS WITH THE RNP.

We will show that if $G$ is second countable and weakly amenable then $\forall 1<p<\infty$, $A_p^t(G)$ having the RNP for $t=s$ implies that it has it for all $1\leq t\leq s$, where $s=\infty$ is allowed.

\medskip

\indent {\it{Definition:}} {\it{Let}} $X$, $Y$ {\it{be Banach spaces and}} $T:X\rightarrow Y$ {\it{be a bounded linear operator.}} $T$ {\it{is a semi-embedding if it is one to one and it maps the closed unit ball in X into a closed set in Y.}} {\it{If such T exists we say that X semi embeds in Y.}}

\medskip

\noindent {\it{THEOREM (H.P. Rosenthal):}} {\it{A separable Banach space has the RNP if it semi-embeds in a Banach space with the RNP.}}  See \cite{Uhl} p.160 or \cite{Ros}, \cite{lpp}.
 
\noindent We will use of the above Theorem, to prove the main Theorem~C.

\noindent We need the following:

\medskip
  
\noindent {\it{LEMMA 4:}} {\it{If}} $r<s$  {\it{then the identity}} $I:A_p^r(G)\rightarrow A_p^s(G)$ {\it{is a semi-embeding, for any}} $s\leq\infty$ {\it{(If}} $s=\infty$, $A_p^\infty (G)=A_p(G)$).

\medskip

\noindent {\it{PROOF:}} Denote by $B_t$ the closed unit ball of $A_p^t(G)$. Let $v_n\in B_r$ satisfy that $\| v_n-w\|_{A_p^s} = \| v_n-w\|_{A_p}+\| v_n-w\|_{L^s}\rightarrow 0$, for some $w\in A_p^s(G)$.

\noindent (If $s=\infty$ only $\| v_n-w\|_{A_p}$ appears). Clearly $|v_n(x)|\rightarrow|w(x)|$, $\forall x\in G$. And by Fatou's Lemma we have $\int |w|^r\, dx\leq\lim\inf\int |v_n|^r\, dx\leq 1$. Thus $w\in A_p^r$. But $1\geq\lim\sup\left(\| v_n\|_{A_p}+\| v_n\|_{L^r}\right)\geq\lim \| v_n\|_{A_p}+\lim\inf \| v_n\|_{L^r}\geq\| w\|_{A_p}+\| w\|_{L^r}$. Thus $w\in B_r$.$\quad\square$ 

\medskip

\noindent {\it{THEOREM C:}} {\it{Assume that G is second countable and weakly amenable. If for some}} $t\leq\infty$, $A_p^t(G)$ has the RNP, then so does $A_p^r(G)$, $\forall 1\leq r\leq t$. {\it{In particular, if}} $A_p(G)$ has the RNP then $A_p^r(G)$ {\it{has the RNP for all}} $1\leq r<\infty$. 

\medskip

\noindent {\it{PROOF:}} Apply Rosenthal's Theorem and the above Lemma~4, and note that by Lemma~1, $A_p^r(G)$ is norm separable, since $\| u\|_\infty \leq \| u\|_{B_p}\leq \| u\|_{B_2}$. $\quad\square$

\medskip

\noindent REMARKS: (1) A group $G$ with completely reducible regular representation is called in \cite{Tlr} an $[AR]$ {\it{group. G is such iff}} $A_2(G)$ {\it{has the RNP}}, as proved by Keith Taylor \cite{Tlr}. A noncompact [AR] group is called a Fell group in \cite{Bag}. Larry Baggett and Keith Taylor construct in \cite{Tay} p.596 (iii) an example of a connected nonunimodular Lie group $G=G_3$ such that $A_2(G)\neq W_2\cap C_0(G)$ and such that $G$ is a Fell group. The above implies that $A_2^r(G)$ {\it{has the RNP, for all}} $1\leq r\leq \infty$ for the above and any Fell group. \cite{Tay} includes examples of Fell groups which are connected Lie groups and which are (i) solvable, (ii) amenable nonsolvable, (iii) nonamenable, (iv) non-TypeI. All of which are not unimodular. See also \cite{Mau}.

\indent (2) Assume that for arbitrary $G$, $A_2^s(G)$ having the RNP for some $s>2$ implies the equality $W_2^s(G)=A_2^s(G)$. {\it{It would then follow for}} $G=Z$, {\it{that}} $A_2^s(Z)$ {\it{does not have the RNP for this s.}}  {\it{This is implied by the fact that}} $A_2^s(Z)\neq W_2^s(Z)$, $\forall s>2$ {\it{as proved in \cite{Hew}. Hence there would be no need to take}} $G=SL(2,R)$ {\it{in the remark above, and}} $Z$ {\it{would suffice.}}

\medskip

\noindent QUESTION: {\it{If}} $G$ {\it{is noncompact abelian then}} $A_2(G)$ {\it{does not have the RNP, since its closed unit ball has no extreme points see}} \cite{Die} p.219. Yet, $A_2^r(G)$, $\forall 1\leq r\leq 2$ {\it{does have the RNP, by Theorem~B.}}

\indent {\it{For such $G$ (or even for}} $G=Z$ {\it{and}} $p=2$), {\it{does}} $A_2^r(G)$ {\it{fail to have the RNP for}} $r>2$?

\bigskip

\noindent (III) THE SCHUR AND THE DUNFORD-PETTIS PROPERTY.

The following result clarifies the DPP case, for discrete $G$ and $1\leq r\leq\max \{ p,p^\prime\}$.

\medskip

\noindent {\it{PROPOSITION 5:}} {\it{Let $G$ be any discrete group. Then}}
\begin{enumerate}
\item {\it{For any}} $1<r\leq\max (p,p^\prime )$, $A_p^r(G)$ {\it{fails the}} DPP {\it{and afortiori fails the SP, yet has the RNP.}}
\item $A_p^1(G)$ {\it{has the SP and the RNP.}}
\end{enumerate}

\medskip

\noindent {\it{PROOF:}} (1) By Theorem~7 of \cite{Gr3} $A_p^r(G)=\ell^r (G)$, if $1\leq r\leq\max\{ p,p^\prime\}$. 

\medskip

\noindent Assume in addition that $r>1$. Then, $\ell^r$, $\ell^{r^\prime}$, considered over the positive integers, are isometric to subspaces of $\ell^r(G)$, $\ell^{r^\prime}(G)$, respectively. Let $x_n=(0,0,\ldots 1,0,0,\ldots )$, where $1$ appears in the $n$-th place, considered as an element of $\ell^r$, and let $x_n^\ast$, be defined exactly as $x_n$, but considered as an element of $\ell^{r^\prime}$. Then $x_n$, $x_n^\ast$ are weakly null sequences in $\ell^r$, $\ell^{r^\prime}$, respectively, yet $(x_n,x_n^\ast )=1$. It follows that $A_p^r(G)=\ell^r (G)$, {\it{if}} $1<r\leq\max (p,p^\prime )$, fails the DPP, afortiori the SP. These have the RNP, since they are reflexive Banach spaces.

\noindent (2) If $r=1$ then $\ell^1 (G)$ has the SP and the RNP.

\medskip

\noindent QUESTION: {\it{It is not clear to us if}} $A_p^r(G)$ {\it{has the DPP if}} $r>\max\{ p,p^\prime\}$ {\it{in case}} $G$ {\it{is discrete, or even if}} $G=Z$ {\it{and}} $p=2$.

\indent If $G$ is abelian {\it{non compact then}} $A_2(G)$ {\it{does not have the SP, since}} $\hat{G}$ {\it{is non discrete.}}
 This is substantially improved in the main result, namely

\medskip

\noindent {\it{THEOREM A:}} {\it{Let G be a noncompact locally compact group.}} 
\begin{enumerate}
\item $\forall 1<r\leq\infty$, $A_p^r(G)$ {\it{does not have the SP}}

\hspace{0.5mm} {\it{Yet}} $A_p^1(G)$ {\it{has the SP if G is discrete.}}
				
\item  {\it{If G is weakly amenable then}} $\forall 1<r\leq p^\prime$, $A_p^r(G)$ {\it{does not have the DPP}}

\hspace{0.5mm} {\it{yet it has the RNP,}}
\end{enumerate}

\medskip

\noindent {\it{CONJECTURE:}} {\it{If G is a connected semisimple Lie group with finite center then}} $A_2^r(G)$ {\it{has the DPP for any}} $r>2$.

\medskip

\noindent For a proof, we need the following results.

\medskip
  
\noindent {\it{LEMMA 6:}} {\it{Let G be a locally compact group, and}} $1<p<\infty$.

{\it{If}} $u_n\in L^{p^\prime}$ {\it{and}} $u_n\rightarrow 0$ {\it{weakly in}} $L^{p^\prime}$, {\it{then}}
\begin{equation}
\forall v\in L^p,\, u_n \ast v^\vee\rightarrow 0\mbox{\, {\it{weakly}}\,}=\sigma (A_p,PM_p)\mbox{\, {\it{in}}}\, A_p.
\end{equation}

\medskip

\noindent {\it{PROOF:}} It is enough to prove that
\begin{equation}
(i)\, \forall \Phi\in PM_p, u\in L^{p^\prime},\, v\in L^p\quad \left(\Phi ,u\ast v^\vee\right) =(\Phi\ast v,u).
\end{equation}

\medskip

\noindent Since then $(\Phi ,u_n\ast v^\vee )=(\Phi\ast v,u_n)\rightarrow 0$, since $\Phi\ast v\in L^p$, $\forall v\in L^p$.  Now (i) is an old result of Eymard \cite{Eymar}, a paper, not easily available. Here is a proof based on \cite{Her1}. 

\noindent Any $\Phi\in PM_p$ is in the ultrastrong closure of $PF_p$ in $PM_p$. Let $\Phi\in PM_p$ and $u\in L^{p^\prime}$, $v\in L^p$. Let $w_\alpha \in L^1$, $w_\alpha\rightarrow\Phi$, ultrastrongly. Thus $\| w_\alpha\ast h-\Phi\ast h\|\rightarrow 0$, $\forall h\in L^p$. But then $(w_\alpha ,s)\rightarrow (\Phi ,s)$, $\forall s\in A_p$, since $w_{\alpha\rightarrow\Phi}$, ultraweakly $=\sigma (PM_p,A_p)$. Hence
\begin{equation}
\left(\Phi , u\ast v^\vee\right)\leftarrow\left( w_\alpha ,u\ast v^\vee\right) 
=(w_\alpha\ast v,u)\rightarrow (\Phi\ast v,u).\quad\square
\end{equation}
Recall that $\ell_\alpha u(x)=u(\alpha x)$ if $u$ is a function on $G$.

\medskip

\noindent {\it{LEMMA 7:}} {\it{(i) Let G be}} $\sigma$ {\it{compact and}} $V=V^{-1}$ {\it{be a neighborhood of}} $e$ {\it{such that}} $\overline{V}$ {\it{is compact. Let}} $g_i\rightarrow\infty$ {\it{and}} $u_i=\ell_{g_i}(1_V\ast 1_V)$. {\it{Then}} $\forall 1<r\leq\infty$, $\{ u_i\}$ {\it{is weakly convergent, but not norm convergent in}} $A_p^r(G)$.

\indent {\it{(ii) If G is any noncompact group, then}} $A_p^r(G)$ {\it{does not have the SP}} $\forall 1<r\leq\infty$.

\medskip

\noindent {\it{PROOF:}} Let $u_i=\ell_{g_i}(1_V\ast 1_V)$. Then $\ell_{g_i}1_V\rightarrow 0$ weakly in $L^{p^\prime}$. Since, if $f\in L^{p^\prime}$ and $\epsilon >0$, let $K\subset G$ be compact such that ${\left(\int\limits_F |f|^{p^\prime}\right)}^{1/p^\prime}< \epsilon$ where $F=G\sim K$. Let $k$ satisfy that if $i>k$ then $g_i^{-1}K\subset F$. Then
\begin{equation}
\left|\int f\ell_{g_i}1_V\right|\leq\epsilon\lambda (V)^{1/p}\mbox{ {\it{if}} }i>k.
\end{equation}
By the above Lemmas $u_i\rightarrow 0$ weakly in $A_p$, and clearly weakly in $L^r$, $\forall 1<r<\infty$. (If $r=\infty$, $A_p^r=A_p$). Thus $u_i\rightarrow 0$, {\it{weakly in}} $A_p^r$, $\forall 1<r\leq\infty$, by Cor.~6 in \cite{Gr3}. But $1_V\ast 1_V(x)=\lambda\{ xV\cap V\} =\lambda\{ V\}$, {\it{if}} $x=1$. Hence
\begin{equation}
{\| u_i\|}_{A_p^r}\geq {\| u_i\|}_\infty =\lambda\{ V\} >0. 
\end{equation}
Thus $\{ u_i\}$ is not norm convergent in $A_p^r(G)$. If $G$ is not $\sigma$ compact, let $H$ be a $\sigma$ compact open subgroup. Consider the two Banach algebras $A=A_p^r(H)$ and $B=\{ 1_Hu;u\in A_p^r(G)\}$ as a closed subalgebra of $A_p^r(G)$. Both these are semisimple, since $A_p(H)$ and $A_p(G)$ are such, see \cite{Her1} Proposition~5. Define $T:B\rightarrow A$, by $T(1_Hv)=v$. Then $T$ is a 1-1 onto algebraic isomorphism. By \cite{Rud} Thm.~11.10, both $T^{-1}$ {\it{and}} $T$ are continuous. Since $H$ is $\sigma$ compact let $\{ u_i\}\subset A$ be a sequence which converges weakly, but not in norm, in $A$. Then $\left\{ T^{-1}u_i\right\}$, satisfies the same in $B$ hence in $A_p^r(G)$. $\quad\square$

\medskip

\noindent {\it{PROOF OF THEOREM A:}} (i) is proved in the above Lemma.

\noindent (ii) If $1\leq r\leq p^\prime$ and $G$ is weakly amenable then $A_p^r(G)$ is a dual Banach space with the RNP. Thus $A_p^r(G)=B^\ast$ for some Banach space $B$, and $B$ does not contain $\ell^1$ by \cite{Ha}. If $A_p^r(G)$ has the DPP then so does $B$, by \cite{Dtl} Cor.~2. But then $B^\ast$ has the SP by \cite{Dtl} Thm.~3, p.23, which contradicts (i), if $r>1$.

\medskip

\noindent {\it{CONJECTURE:}} {\it{If G is a connected semisimple real Lie group then}} $A_2^r(G)$ {\it{for}} $r>2$, {\it{has the DPP.}}

\bigskip

\centerline{Department of Mathematics, University of British Columbia}

\centerline{Vancouver B.C. V6T 1Z4, Canada}

\centerline{E-mail address: granirer@math.ubc.ca}
\end{document}